\documentclass[a4paper,10pt]{article}
\usepackage{amsmath}
\usepackage{amssymb}
\usepackage{amsthm}
\theoremstyle{plain}
\numberwithin{equation}{section}
\newtheorem{prop}{Proposition}[section]
\newtheorem{theo}{Theorem}[section]
\newtheorem{cor}{Corollary}[section]
\newtheorem{lem}{Lemma}[section]
\theoremstyle{definition}
\newtheorem{defi}{Definition}[section]
\theoremstyle{remark}
\newtheorem{rem}{Remark}[section]

\begin{document}
\begin{Large}
 {\bf Subexponential densities of infinitely 
divisible\\
 distributions on the half line}
\end{Large}
\begin{center}
  Toshiro Watanabe\\
Center for Mathematical Sciences, The University of Aizu,\\
Ikkimachi Tsuruga, Aizu-Wakamatsu, Fukushima 965-8580, Japan \\
E-mail: markov2000t@yahoo.co.jp
\end{center}  
\medskip
{\bf Abstract.}  We show that, under the long-tailedness of the densities of normalized L\'evy measures, the densities  of infinitely divisible distributions on  the half line  are subexponential if and only if the densities  of their normalized L\'evy measures are  subexponential. Moreover, we prove that, under a certain continuity assumption, the densities  of infinitely divisible distributions on  the half line  are subexponential if and only if  their normalized L\'evy measures are locally subexponential. 
  
\medskip

  infinite divisibility, L\'evy measure\\
 

{\bf Abbreviated title} : Subexponential densities
 
\section{Introduction and main results}
\medskip
The subexponentiality of infinitely divisible distributions on the half line was characterized by Embrechts et al. \cite{egv} and on the real line by Pakes \cite{p} and Watanabe \cite{w}.  The subexponentiality of a density is a stronger and  more difficult property than the subexponentiality of a distribution. Some infinitely divisible distributions  on the half line  such as Pareto, lognormal, and Weibull (with parameter less than 1) distributions have subexponential densities. Watanabe and Yamamuro \cite{wy2} proved that the density of a self-decomposable distribution on the real line is subexponential if and only if the  density on $(1,\infty)$ of its  normalized  L\'evy measure is subexponential. The purpose of this paper is to characterize the subexponential densities of  absolutely continuous infinitely divisible distributions on the half line under some additional assumptions.  
 
  In what follows, we denote by $\mathbb R$ the real line and by $\mathbb R_{+}$ the half line $[0,\infty)$. Denote by $\mathbb N$ the totality of positive integers. The symbol $\delta_a(dx)$ stands for the delta measure at $a \in \mathbb R$. Let $\eta$ and $\rho$ be finite measures on $\mathbb R$. We denote by $\eta*\rho$ the convolution of $\eta$ and $\rho$ and by $\rho^{n*}$
$n$-th convolution power of $\rho$ with the understanding that $\rho^{0*}(dx)=\delta_0(dx)$. The characteristic function of a distribution $\rho$ is denoted by $\widehat\rho(z)$, namely, for $z \in \mathbb R$, 
\begin{equation}
\widehat\rho(z)=\int_{-\infty}^{\infty}e^{izx}\rho(dx). \nonumber
\end{equation}
Let $f(x)$ and $g(x)$ be probability density functions on $\mathbb R$. We denote by $f \otimes g(x)$ the convolution of $f(x)$ and $g(x)$ and by $f^{n\otimes}(x)$ $n$-th convolution power of $f(x)$ for $n \in \mathbb N$. For positive functions $f_1(x)$ and $g_1(x)$ on $[A,\infty)$ for some $A \in \mathbb R$, we define the relation $f_1(x)\sim g_1(x)$ by $\lim_{x \to \infty}f_1(x)/g_1(x)=1.$   We use the symbols $\mathcal{L}$ and $\mathcal{S}$ in the sense of long-tailed and subexponential, respectively.
\begin{defi}
(i) A nonnegative measurable function $g(x)$ on $\mathbb R$ belongs to the class ${\bf L}$ if $g(x+a)\sim g(x)$ for every $ a \in \mathbb R$.

(ii) A probability density function $g(x)$ on $\mathbb R$ belongs to the class $\mathcal{L}_d$ if $g(x) \in {\bf L}$. A distribution $\rho$ on $\mathbb R$ belongs to the class $\mathcal{L}_{ac}$ if there is $g(x) \in \mathcal{L}_d$ such that $\rho(dx)=g(x)dx$.

(iii)  A probability density function $g(x)$ on $\mathbb R$ belongs to the class $\mathcal{S}_d$ if $g(x) \in \mathcal{L}_d$ and $g^{2\otimes}(x)\sim 2g(x)$.  A distribution $\rho$ on $\mathbb R$ belongs to the class $\mathcal{S}_{ac}$ if there is $g(x) \in 
\mathcal{S}_d$ such that $\rho(dx)=g(x)dx$.
\end{defi}
\begin{defi} (i)  Let $\Delta:=(0,c]$ with $c >0$. A distribution $\rho$ on $\mathbb R$ belongs to the class $\mathcal{L}_{\Delta}$ if $\rho((x,x+c]) \in {\bf L}$. A distribution $\rho$ on $\mathbb R$ belongs to the class $\mathcal{L}_{loc}$ if $\rho \in \mathcal{L}_{\Delta}$ for each $\Delta:=(0,c]$ with $c >0$.

(ii)  Let $\Delta:=(0,c]$ with $c >0$. A distribution $\rho$ on $\mathbb R$ belongs to the class $\mathcal{S}_{\Delta}$ if $\rho \in \mathcal{L}_{\Delta}$ and  $\rho^{2*}((x,x+c])\sim 2\rho((x,x+c])$. 
 A distribution $\rho$ on $\mathbb R$ belongs to the class $\mathcal{S}_{loc}$ if $\rho \in \mathcal{S}_{\Delta}$ for each $\Delta:=(0,c]$ with $c >0$.

\end{defi}
\medskip
 A probability distribution $\rho$ on $\mathbb R$ is called {\it subexponential} if $\rho((x,\infty)) \in {\bf L}$ and  $\rho^{2*}((x,\infty))\sim 2\rho((x,\infty))$. The class of all subexponential 
 distributions on $\mathbb R$ is denoted by $\mathcal{S}$. Functions in the class ${\bf L}$ are called {\it long-tailed functions}. Probability density functions in the classes  $\mathcal{L}_d$ and  $\mathcal{S}_d$ are called {\it long-tailed densities} and {\it subexponential densities}, respectively. The class $\mathcal{S}_d$ was introduced by Chover et al. \cite{cnw}.  Note that if $f(x) \in \mathcal{L}_d$, then $\lim_{x \to \infty}f(x)=0$ and $\lim_{x \to \infty}e^{sx}f(x)=\infty$ for every $s >0$. See Foss et al. \cite{fkz}.  Distributions in the classes  $\mathcal{S}_{\Delta}$ and $\mathcal{S}_{loc}$ are called $\Delta$-{\it subexponential} and {\it locally subexponential}, respectively. The class  $\mathcal{S}_{\Delta}$ was introduced by Asmussen et al. \cite{afk}. The class $\mathcal{S}$ includes the classes $\mathcal{S}_{\Delta}$ and $\mathcal{S}_{ac}$. 

 Let $\mu$ be an infinitely divisible distribution on $\mathbb R_+$. Then, its characteristic function $\widehat\mu(z)$ is represented as
\begin{equation}
\widehat\mu(z)=\exp \left( \int_0^{\infty}(e^{izx}-1)\nu(dx)+i\gamma_0z \right),         
\nonumber
\end{equation}
where $\gamma_0 \in \mathbb R_+$ and $\nu$ is a measure on $\mathbb R_+$ satisfying $\nu(\{0\})=0$ and 
$$\int_0^{\infty}(1\wedge x)\nu(dx) < \infty.$$
The measure $\nu$ is called L\'evy measure of $\mu$. Denote by $\mu^{t*}$ $t$-th convolution power of $\mu$ for $t >0$. Then, $\mu^{t*}$ is a distribution of a certain L\'evy process on $\mathbb R_+$. The distribution $\mu$ is called {\it non-trivial} if it is not a delta measure on  $\mathbb R_+$. A non-trivial infinitely divisible distribution $\mu$ on $\mathbb R_+$ is called {\it semi-stable} if, for some $a > 1$, there exist $b > 1$ and $c \in \mathbb R$ such that, for $z \in \mathbb R$, 
\begin{equation}
 \widehat\mu(z)^a=\widehat\mu(bz)e^{icz}. \nonumber
 \end{equation}
 For a semi-stable distribution on  $\mathbb R_+$, the span $b$ and the index $\alpha:=\log a/\log b$ $(0 < \alpha <1)$ are important. A span is not unique and the index is independent of the choice of span. See Sato \cite{s}.  Shimura and Watanabe \cite{sw} proved that a semi-stable distribution on  $\mathbb R_+$ is subexponential if and only if its L\'evy measure is continuous. Watanabe and Yamamuro \cite{wy3}
 studied the tail behaviors of semi-stable distributions on $\mathbb R$.
 
 Let $\mu$ be an infinitely divisible distribution on $\mathbb R_+$ with L\'evy measure $\nu$. Throughout the paper, we assume that $\gamma_0=0$ and that the tail $\nu((c,\infty)) $ is positive for all $c > 0$. Define a normalized distribution $\nu_{(1)}$ as
 $$ \nu_{(1)}(dx):=1_{(1,\infty)}(x)\frac{\nu(dx)}{\nu((1,\infty))}.$$
  Here the symbol $1_{(1,\infty)}(x)$ stands for the indicator function of the set $(1,\infty)$. It is known that if 
  \begin{equation}
\int_{-\infty}^{\infty}|\widehat\mu(z)|dz < \infty,
\end{equation}
then $\mu$ is absolutely continuous with a bounded and continuous density on $\mathbb R$. It is also  known that if $\mu$ is an  infinitely divisible distribution on $\mathbb R_+$  with absolutely continuous L\'evy measure $\nu(dx)$ satisfying $\nu((0,1))=\infty$, then $\mu$ is absolutely continuous. First, we have the following basic result.
 \begin{prop} Let $\mu$ be an absolutely continuous  infinitely divisible distribution on $\mathbb R_+$  with L\'evy measure $\nu$  and let $p(x)$ be a density of $\mu$. Then, we have the following.
 
(i) Assume that $p(x) \in \mathcal{L}_d$. Then, $p(x) \in \mathcal{S}_d$ if and only if  $\nu_{(1)} \in \mathcal{S}_{loc}$.

 (ii)  The following are equivalent : 

(a) $p(x) \in \mathcal{S}_d$.

(b) $\nu_{(1)} \in \mathcal{L}_{loc}$
and $p(x) \sim \nu((x,x+1])$.

(c) $\nu_{(1)} \in \mathcal{L}_{loc}$
and there is $C \in (0,\infty)$ such that $p(x) \sim C\nu((x,x+1])$.
\end{prop}
\medskip
Next, we characterize the case where $\nu_{(1)}$ is absolutely continuous with a long-tailed density. Note that Corollary 1.1 of Jiang et al. \cite{jwcc} is analogous to  Theorem 1.1 below  but their proof is not valid.
\medskip

 \begin{theo}
Let $\mu$ be an  infinitely divisible distribution on $\mathbb R_+$  with absolutely continuous L\'evy measure $\nu(dx)=g(x)dx$ satisfying $\nu((0,1))=\infty$. Let $g_1(x)$ be a density of $\nu_{(1)}$. Let $f_1(x)$ be a density  of the infinitely divisible distribution on $\mathbb R_+$ with L\'evy measure $1_{(0,1)}(x)\nu(dx)$. Assume that for some $b_0 >0$
 $$\lim_{x \to \infty} \exp(b_0x)f_1(x)=0$$
  and that $g_1(x)$ is bounded and $g_1(x) \in \mathcal{L}_{d}$.    Then, we can choose   a density $p(x)$ of $\mu$ such that the following are equivalent : 

(1) $p(x) \in \mathcal{S}_d$.

(2) $g_1(x) \in \mathcal{S}_{d}$.

(3)  $p(x) \sim C_0 g_1(x)$ with $C_0 :=\int_1^{\infty}g(x)dx$.

(4) There is $C \in (0,\infty)$ such that $p(x) \sim Cg_1(x)$.
\end{theo}
\begin{rem} Assume that $\nu(dx) =1_{(0,1)}(x)x^{-1}k(x)dx +1_{(0,\infty)}(x)h(x)dx$, where $k(x)$ is nonnegative and decreasing with $0 <k(0+) \le \infty$ and $\int_0^1k(x)dx < \infty$ and $h(x)$ is nonnegative, bounded, and  integrable on $\mathbb R_+$. Then we can choose   a density  $f_1(x)$ of the infinitely divisible distribution on $\mathbb R_+$ with L\'evy measure $1_{(0,1)}(x)\nu(dx)$ such that, for any $b >0$,
 $$\lim_{x \to \infty} \exp(bx)f_1(x)=0.$$
\end{rem}
\begin{rem} Theorem 1.1 includes all self-decomposable cases  on $\mathbb R_+$ of Watanabe and Yamamuro \cite{wy2}. We see from Theorem 53.6 of Sato \cite{s} that  there is a self-decomposable case  on $\mathbb R_+$ of Watanabe and Yamamuro \cite{wy2}, which does not satisfies the condition  (1.1). Hence Theorem 1.2 below does not include Theorem 1.1.
\end{rem}

\medskip
We 
apply Theorem 1.1 to the
 tail asymptotic behavior of the density of the distribution of a L\'evy process  on $\mathbb R_+$.

\begin{cor} 
Let $\mu$ be an  infinitely divisible distribution on $\mathbb R_+$ with absolutely continuous L\'evy measure $\nu(dx)=g(x)dx$ satisfying $\nu((0,1))=\infty$. . Let $g_1(x)$ be a density of $\nu_{(1)}$.  Let $f_1^t(x)$ be a density  of the infinitely divisible distribution on $\mathbb R_+$ with L\'evy measure $t1_{(0,1)}(x)\nu(dx)$. Assume that, for every $t >0$, there is some $b_0(t) >0$ such that
 $$\lim_{x \to \infty} \exp(b_0(t)x)f_1^t(x)=0$$ 
 and that $g_1(x)$ is bounded and $g_1(x) \in \mathcal{L}_{d}$.   
 Then, we can choose   a density $p^t(x)$ of $\mu^{t*}$ such that the following hold : 
 
(i) If $p^t(x) \in \mathcal{S}_d$ for some $t >0$, then $p^t(x) \in \mathcal{S}_d$ for all $t >0$ and
$$ p^t(x) \sim t p^1(x) \mbox{ for all }t >0. $$

(ii) If $p^1(x) \in \mathcal{L}_d$ and, for some 
$t \in  (0,1) \cup (1,\infty)$, there is 
$ C(t) \in (0,\infty)$            
 such that

\begin{equation}
p^t(x)\sim  C(t)p^1(x),
\end{equation}
then $C(t)=t$ and $p^1(x) \in \mathcal{S}_d$.
\end{cor}

Under the assumption (1.1), we can characterize the subexponential density of an infinitely divisible distribution on $\mathbb R_+$ as follows.

\begin{theo}
Let $\mu$ be an  infinitely divisible distribution on $\mathbb R_+$  with L\'evy measure $\nu$. Assume that (1.1) holds.  Let $p(x)$ be the continuous density of $\mu$. Then, the following are equivalent : 

(1) $p(x) \in \mathcal{S}_d$.

(2) $\nu_{(1)} \in \mathcal{S}_{loc}$.

(3) $\nu_{(1)} \in \mathcal{L}_{loc}$
and $p(x) \sim \nu((x,x+1])$.

(4) $\nu_{(1)} \in \mathcal{L}_{loc}$
and there is $C \in (0,\infty)$ such that $p(x) \sim C\nu((x,x+1])$.

\end{theo}
\begin{rem} (i) Assume that $\nu(dx)=g(x)dx$ with $\liminf_{x \to 0+}xg(x) >1$, then  we have (1.1). See Lemma 53.9 of Sato \cite{s}. Some sufficient conditions in order that certain  infinitely divisible distributions $\mu$  on $\mathbb {R}_+$ satisfy (1.1) are found in Sato \cite{s} and Watanabe \cite{w3}. 

(ii) Let $\mu$ be a semi-stable distribution on $\mathbb R_+$ with  L\'evy measure $\nu$. By virtue of Proposition 24.20 of Sato \cite{s}, we have (1.1). Thus, we see from the above theorem that $\mu \in   \mathcal{S}_{ac}$ if and only if $\nu_{(1)} \in \mathcal{S}_{loc}$.

(iii) Let $1 < x_0 < b$ and choose $\delta \in (0,1)$ satisfying $\delta < (x_0 -1)\wedge (b-x_0).$ We take a continuous periodic function $h(x)$ on $\mathbb R$ with period $\log b$ such that $h(\log x)>0$ for $x \in [1,x_ 0)\cup (x_0,b]$ and
\begin{equation}
h(\log x):=\left \{
\begin{array}{l}
0\quad \quad \quad \ \mbox{ for } x=x_0,\\
\frac{-1}{\log |x-x_0|}\mbox{ for } 0 < |x-x_0|< \delta.
\end{array}
\right. \nonumber
\end{equation}
Let $\Phi(x):=x^{-\alpha-1}h(\log x)$ on $(0,\infty)$ with $\alpha \in (0,1).$ Then $\nu$ defined by
$$\nu(dx):=1_{(0,\infty)}(x)\Phi(x)dx$$ 
is the L\'evy measure of a certain semi-stable distribution on $\mathbb R_+$ with span $b$ and index $\alpha$.
Watanabe and Yamamuro proved in the proof of Theorem 1.1 of \cite{wy4} that $\nu_{(1)} \in \mathcal{S}_{loc}\setminus \mathcal{S}_{ac}$. Thus, we find from (ii) that $\mu \in \mathcal{S}_{ac}$ but $\nu_{(1)} \notin    
 \mathcal{S}_{ac}$ in this case and hence Theorem 1.1 does not include Theorem 1.2. Note that the continuous density $p(x)$ of this semi-stable distribution $\mu$ is not almost decreasing because
$$
\lim_{n \to \infty}\frac{p(b^n(x_0+\delta))}{p(b^nx_0)}=\lim_{n \to \infty}\frac{\nu((b^n(x_0+\delta),b^n(x_0+\delta)+1])}{\nu((b^nx_0,b^nx_0+1])}=\infty.$$
\end{rem}
\medskip

Finally, we can
apply the theorem above to the
 tail asymptotic behavior of the density of the distribution of a L\'evy process  on $\mathbb R_+$.
\begin{cor} 
Let $\mu$ be an  infinitely divisible distribution on $\mathbb R_+$.  
 Assume that $\int_{-\infty}^{\infty}|\widehat\mu(z)|^tdz < \infty$ for every $t >0$.  Let $p^t(x)$ be the continuous density of $\mu^{t*}$ for $t >0.$  Then, we have the following.

(i) If $p^t(x) \in \mathcal{S}_d$ for some $t >0$, then $p^t(x) \in \mathcal{S}_d$ for all $t >0$ and
$$ p^t(x) \sim t p^1(x) \mbox{ for all }t >0. $$

(ii) If $p^1(x) \in \mathcal{L}_d$ and, for some 
$t \in  (0,1) \cup (1,\infty)$, there is 
$ C(t) \in (0,\infty)$            
 such that

\begin{equation}
p^t(x)\sim  C(t)p^1(x),
\end{equation}
then $C(t)=t$ and $p^1(x) \in \mathcal{S}_d$.
\end{cor}

The organization of this paper is as follows. In Sect.\ 2, we explain basic results on the classes $\mathcal{L}_d$,  $\mathcal{S}_d$,
and $\mathcal{S}_{loc}$ as preliminaries. In Sect.\ 3,  we prove the main results.

\section{Preliminaries}

 In this section, we give several fundamental results on the classes $\mathcal{L}_d$, $\mathcal{S}_d$,
and $\mathcal{S}_{loc}$.

\begin{lem} Let $f(x)$ and $g(x)$ be probability density functions on $\mathbb R$.

(i) If $f(x),g(x) \in \mathcal{L}_d$, then $f\otimes g(x) \in \mathcal{L}_d$.

(ii) Let $f(x) \in \mathcal{L}_d$ and define a distribution $\rho$ on $\mathbb R$ by 
$$\rho(dx):=f(x)dx.$$ 
 Then, $\rho \in \mathcal{S}_{loc}$ if and only if $f(x) \in \mathcal{S}_d$. That is, $\mathcal{S}_{ac}= \mathcal{S}_{loc}
\cap \mathcal{L}_{ac}.$
\end{lem}
\medskip

Proof.   Assertion (i) is due to Theorem 4.3 of Foss et al. \cite{fkz}. Next, we prove assertion (ii). Assume that $f(x) \in \mathcal{L}_d$. Then, by (i), $f^{2\otimes}(x) \in \mathcal{L}_d$
and hence $f(x+u) \sim f(x)$ and $f^{2\otimes}(x+u) \sim f^{2\otimes}(x)$ uniformly in $u \in [0,c]$ with $c >0$. Thus, we have, for $x > 0$,
$$\rho((x,x+c])=\int_0^{c}f(x+u)du\sim cf(x)$$
and
\begin{equation}
\rho^{2*}((x,x+c])=\int_0^{c}f^{2\otimes}(x+u)du\sim cf^{2\otimes}(x). \nonumber
\end{equation}
Hence, we see that  $\rho \in \mathcal{S}_{loc}$ if and only if $f^{2\otimes}(x) \sim 2f(x)$, namely, $f(x) \in \mathcal{S}_d$. 
$\Box$

\begin{lem}  Let $f(x)$ and $g(x)$ be probability density functions on $\mathbb R_+$.

(i) If $f(x) \in \mathcal{S}_d$ and $g(x) \sim cf(x)$ with $c \in (0,\infty)$, then $g(x) \in \mathcal{S}_d$.

(ii) Assume that $f(x) \in \mathcal{S}_d$ and $f(x)$ is bounded on $\mathbb R_+$. Then, for any $\epsilon > 0$, there are $x_0(\epsilon) >0$ and $C(\epsilon) >0$ such that, for all $x > x_0(\epsilon)$ and all $n \in \mathbb N$, 
$$f^{n\otimes}(x) \leq C(\epsilon)(1+\epsilon)^nf(x).$$

(iii) If $f(x) \in \mathcal{S}_d$, then, for all $n \in \mathbb N$,
$$ f^{n\otimes}(x) \sim nf(x).$$
\end{lem}

\medskip
Proof. Assertions (i), (ii), and (iii) are due to Theorem 4.8, Theorem 4.11, and Corollary 4.10 of Foss et al. \cite{fkz}, respectively.  \hfill $\Box$

\medskip

\begin{lem}
 (i) If a distribution $\rho$ on $\mathbb R_+$ belongs to 
 $\mathcal{L}_{loc}$, then $\rho((x,x+c])\sim c\rho((x,x+1])$ for all $c >0$.
 
 (ii) Let $q(x)$ be a continuous probability density function 
 on $\mathbb R_+$ such that, for some $\gamma_1 >0$, $\lim_{x \to \infty}e^{\gamma_1 x}q(x)=0.$ Let $\rho$ be a distribution  on $\mathbb R_+$ and define a probability density function $p(x)$ on $\mathbb R_+$ as
 \begin{equation}
 p(x):=\int_{0-}^{x+}q(x-u)\rho(du).
 \end{equation}
Then, $\rho \in \mathcal{L}_{loc}$ implies $p(x) \in \mathcal{L}_{d}$.
 \end{lem}
 \medskip
 
Proof.  Assertion (i) is proved as (2.6) in Theorem 2.1 of Watanabe and Yamamuro \cite{wy2}. Next, we prove assertion (ii). Suppose that $\rho \in \mathcal{L}_{loc}$ and a probability density function $q(x)$ on $\mathbb R_+$ is continuous on $\mathbb R_+$ and  $\lim_{x \to \infty}e^{2\gamma x}q(x)=0$ for some $\gamma >0$. Let  $p(x)$ be  a probability density function  on $\mathbb R_+$ defined by (2.1).  Let $N \in \mathbb N$ and $x > 2N. $
 We have $p(x)=\sum_{j=1}^3I_j(x)$, where
 $$ I_1(x):=\int_{(x-N)+}^{x+}q(x-u)\rho(du),$$
  $$ I_2(x):=\int_{0-}^{N+}q(x-u)\rho(du),$$
 and
  $$ I_3(x):=\int_{N+}^{(x-N)+}q(x-u)\rho(du).$$
 For $M \in \mathbb N$, there are $\delta(M)\geq 0$ and $a_n=a_n(M)\geq 0$ for all $n \in \mathbb N$ such that 
 $$a_n \leq q(x) \leq a_n+ \delta(M)$$
 for $M^{-1}(n-1)\leq x \leq M^{-1}n$ and $1\leq n \leq MN$. Define $J(x;M,N)$ as
 $$J(x;M,N):=\sum_{n=1}^{MN}a_n\rho((x-M^{-1}n,x-M^{-1}(n-1)]).$$
 Then, we see from (i) that
 $$J(x;M,N)\sim \rho((x, x+1])\sum_{n=1}^{MN}a_nM^{-1}$$
 and
 $$J(x;M,N) \leq I_1(x) \leq J(x;M,N)+ \delta(M)\rho((x-N,x]).$$
 Since we can make $\lim_{M \to \infty}\delta(M)=0$ and
 $$ \lim_{M \to \infty}\sum_{n=1}^{MN}a_nM^{-1}=\int_0^Nq(x)dx,$$
 we find that
 \begin{equation}
 I_1(x)\sim\rho((x, x+1])\int_0^Nq(x)dx.
 \end{equation}
 For $x \in \mathbb R$, the symbol $[x]$ stands for the largest integer not exceeding $x$. We can take $N \in \mathbb N$ sufficiently large such that
 $q(y) \leq e^{-2\gamma y}$ for $y > N$ and, for some $C_1>0$ and for integers $n$ in $0 \leq n \leq [x]-2N+1$ with $x > 2N$,
 $$ \rho((N+n,N+n+1]) \leq C_1e^{\gamma(x-N-n)}\rho((x,x+1]).$$
 Thus, we have
 \begin{equation}
 \begin{split}
 I_2(x) &\leq \int_{0-}^{N+}e^{-2\gamma(x-u)}\rho(du) \\ 
 &\leq e^{-2\gamma(x-N)}=o(\rho((x, x+1])).
 \end{split}
 \end{equation}
Moreover, we obtain that, for some $C_2 >0$, 
  \begin{equation}
 \begin{split}
 I_3(x) &\leq \int_{N+}^{(x-N)+}e^{-2\gamma(x-u)}\rho(du) \\ \nonumber
 &\leq\sum_{n=0}^{[x]-2N+1} e^{-2\gamma(x-N-n-1)}\rho((N+n,N+n+1]))
 \\
 &\leq e^{2\gamma}C_1\sum_{n=0}^{[x]-2N+1}e^{-\gamma(x-N-n)}\rho((x,x+1])\\
&\leq C_2e^{-\gamma N}\rho((x,x+1]) .
\end{split}
\end{equation}
Thus, we see that
\begin{equation}
\lim_{N \to \infty}\limsup_{x \to \infty}\frac{I_3(x)}{\rho((x,x+1])}=0.
\end{equation}
Hence, we find from (2.2), (2.3), and (2.4) that
$$ p(x) \sim \rho((x,x+1])\int_0^{\infty}q(x)dx=\rho((x,x+1])$$
and thereby $p(x) \in \mathcal{L}_{d}$.\hfill $\Box$
%\end{proof}
\medskip

Watanabe and Yamamuro \cite{wy2} used the  results of Watanabe \cite{w} on the convolution equivalence of infinitely divisible distributions on $\mathbb R$  to prove the following  lemmas. Our main results essentially depend on those two  results.

\begin{lem} (Theorem 1.1 of \cite{wy2}) Let $\mu$ be an  infinitely divisible distribution on $\mathbb R_+$  with L\'evy measure $\nu$. Then, the following are equivalent : 

(1) $\mu \in \mathcal{S}_{loc}$.

(2) $\nu_{(1)} \in \mathcal{S}_{loc}$.

(3) $\nu_{(1)} \in \mathcal{L}_{loc}$
and $\mu((x,x+c]) \sim \nu((x,x+c])$ for all $c >0$.

(4) $\nu_{(1)} \in \mathcal{L}_{loc}$
and there is $C \in (0,\infty)$ such that $\mu((x,x+c]) \sim C\nu((x,x+c])$ for all $c >0$.
\end{lem}

\begin{lem} (Theorem 1.2 of \cite{wy2}) Let $\mu$ be an  infinitely divisible distribution on $\mathbb R_+$  with L\'evy measure $\nu$. Then, we have the following.

(i) If $\mu^{t*} \in \mathcal{S}_{loc}$ for some $t >0$, then $\mu^{t*} \in \mathcal{S}_{loc}$ for all $t >0$ and
$$ \mu^{t*}((x,x+c]) \sim t\mu((x,x+c]) $$
 for all $t >0 $ and for all $c >0. $

(ii) If $\mu \in \mathcal{S}_{loc}$ and, for some 
$t \in  (0,1) \cup (1,\infty)$, there is 
$ C(t) \in (0,\infty)$            
 such that

\begin{equation}
 \mu^{t*}((x,x+c]) \sim C(t)\mu((x,x+c])
\end{equation}
 for all $c >0, $ then $C(t)=t$ and $\mu \in \mathcal{S}_{loc}$.
\end{lem}

\section{Proofs of the main results}

In this section, let $\mu$ be an  infinitely divisible distribution on $\mathbb R_+$  with L\'evy measure $\nu$ and let $\mu_1$ be a compound Poisson  distribution on $\mathbb R_+$  with L\'evy measure $\nu_1(dx):=1_{(c_1,\infty)}(x)\nu(dx)$ for $c_1 >0$. We define  an  infinitely divisible distribution $\zeta_1$ on $\mathbb R_+$ by $\mu =\mu_1*\zeta_1$.

\begin{lem} (i) We have $\int_{-\infty}^{\infty}|\widehat\mu(z)|dz < \infty$ if and only if $\int_{-\infty}^{\infty}|\widehat\zeta_1(z)|dz < \infty$ for all sufficiently large $c_1 >0$.

(ii) Suppose that $\int_{-\infty}^{\infty}|\widehat\zeta_1(z)|dz < \infty$ for some $c_1 >0$. Let $q(x)$ be the continuous density of $\zeta_1$ on $\mathbb R_+$. Then, for any $\gamma > 0$, $\lim_{x \to \infty}e^{\gamma x}q(x)=0.$
\end{lem}
\medskip

Proof. First, we prove assertion (i). Since 
$$| \widehat\mu(z)|=|\widehat\mu_1(z)||\widehat\zeta_1(z)|$$
 and $|\widehat\mu_1(z)|\leq 1$, we have $| \widehat\mu(z)|\leq |\widehat\zeta_1(z)|$. Thus, if $\int_{-\infty}^{\infty}|\widehat\zeta_1(z)|dz < \infty$ for some $c_1 >0$, then $\int_{-\infty}^{\infty}|\widehat\mu(z)|dz < \infty$. Let
 $\lambda:=\nu((c_1,\infty))$ and define a constant $C_{\nu}$ as
 $$C_{\nu}:=\inf\{c >0:\nu((c,\infty)) < \log 2 \}.$$
Since $\mu_1$ is a compound Poisson distribution and 
$$|\widehat\mu_1(z)-e^{-\lambda}| = |\widehat\mu_1(z)-e^{-\lambda} \widehat\delta_0(z)|\le 1-e^{-\lambda},$$
  we see that, for all $c_1 > C_{\nu}$,
 $$|\widehat\mu_1(z)|\geq e^{-\lambda}-|\widehat\mu_1(z)-e^{-\lambda}|\geq 2e^{-\lambda}-1 >0.$$
 Hence we have $| \widehat\mu(z)|\geq C_0|\widehat\zeta_1(z)|$ with some $C_0>0$ for all $c_1 > C_{\nu}$. Thus, if $\int_{-\infty}^{\infty}|\widehat\mu(z)|dz < \infty$, then $\int_{-\infty}^{\infty}|\widehat\zeta_1(z)|dz < \infty$  for all $c_1 > C_{\nu}$. Next, we prove assertion (ii). Suppose that $\int_{-\infty}^{\infty}|\widehat\zeta_1(z)|dz < \infty$ for some $c_1 >0$. It is clear that $q(x)$ is bounded and continuous on $\mathbb R$. Since, for all $\gamma >0$, 
 $$ \int_{1-}^{(c_1\vee 1)+}e^{\gamma x}\nu(dx)< \infty,$$
we see from Theorem 25.3 of Sato \cite{s} that, for all $\gamma >0$,   
  $$ \int_{0-}^{\infty}e^{\gamma x}\zeta_1(dx)< \infty.$$
Define an exponential tilt $\xi$ of the distribution $\zeta_1$ on $\mathbb R_+$ as
$$\xi(dx):=\frac{e^{\gamma x}q(x)}{ \int_{0-}^{\infty}e^{\gamma x}\zeta_1(dx)}dx.$$
We find from Example 33.15 of Sato \cite{s} that $\xi$ is an  infinitely divisible distribution on $\mathbb R_+$  with L\'evy measure $1_{(0,c_1]}(x)e^{\gamma x}\nu(dx).$ Thus, we have
$$|\widehat\xi(z)|=|\widehat\zeta_1(z)|\exp\left(\int_{0}^{c_1+}(\cos zx -1)(e^{\gamma x} -1)\nu(dx)\right)\leq |\widehat\zeta_1(z)|.$$
 Hence $\int_{-\infty}^{\infty}|\widehat\xi(z)|dz < \infty$ and
  by the Riemann-Lebesgue lemma,  we have, for any $\gamma > 0$, $\lim_{x \to \infty}e^{\gamma x}q(x)=0.$\hfill $\Box$

\medskip
Proof of Proposition 1.1.   We see from Lemma 2.4 that the four conditions (1) $\mu \in \mathcal{S}_{loc}$; (2) $\nu_{(1)} \in \mathcal{S}_{loc}$; (3) $\nu_{(1)} \in \mathcal{L}_{loc}$
and $\mu((x,x+c]) \sim \nu((x,x+c])$ for all $c >0$; (4) $\nu_{(1)} \in \mathcal{L}_{loc}$
and there is $C \in (0,\infty)$ such that $\mu((x,x+c]) \sim C\nu((x,x+c])$ for all $c >0$ are equivalent. First, we prove assertion (i). Assume that $p(x) \in \mathcal{L}_d$. If $p(x) \in \mathcal{S}_d$, then   $\mu \in \mathcal{S}_{loc}$ and hence $\nu_{(1)} \in \mathcal{S}_{loc}$. Conversely, if  $\nu_{(1)} \in \mathcal{S}_{loc}$, then  $\mu \in \mathcal{S}_{loc}$ and thereby $p(x) \in \mathcal{S}_d$ by (ii) of Lemma 2.1. Next, we prove assertion (ii). We obtain from (ii) of Lemma 2.1 that (a), (b), and (c) are respectively equivalent to (1) and $p(x) \in \mathcal{L}_d$; (3) and $p(x) \in \mathcal{L}_d$; (4) and $p(x) \in \mathcal{L}_d$. Hence, 
we find that (a)-(c) are equivalent.\hfill $\Box$

\medskip
Proof of Theorem 1.1.  Assume that for some $b_0 >0$
 $$\lim_{x \to \infty} \exp(b_0x)f_1(x)=0$$
  and that $g_1(x)$ is bounded and $g_1(x) \in \mathcal{L}_{d}$. We see from (i) of Proposition 1.1 and (ii) of Lemma 2.1 that (1) implies (2). We obtain from (ii) of Proposition 1.1 and $g_1(x) \in \mathcal{L}_{d}$ that (1), (3), and (4) are equivalent. Finally, we prove that (2) implies (1).
Let $\nu(dx)=\eta_1(dx)+\eta_2(dx)$, where
$\eta_1(dx):=1_{(0,1)}(x)g(x)dx$ and 
$\eta_2(dx):=1_{(1,\infty)}g(x)dx$. 
 Let $\xi_j(dx)$ be the infinitely divisible distribution on $\mathbb R_+$ with L\'evy measure $\eta_j(dx)$ for $j=1,2$. Then
$\xi_1(dx)=f_1(x)dx$.
By virtue of Theorem 25.3 of Sato \cite{s}, we have, for any $b>0$, 
$$\int_{0-}^{\infty}\exp(bx)f_1(x)dx<\infty.$$
Note that  $\nu_{(1)}(dx)=\eta_2(dx)/\eta_2((1,\infty))$ has the bounded density $g_1(x)$ belonging to  $S_d$ and $\xi_2(dx)$ is a compound Poisson distribution. Thus, we have 
$$\xi_2(dx)=
c\delta_0(dx)+ (1-c)f_2(x)dx,$$
where $c :=\exp(-\eta_2((1,\infty)))\in (0,1)$. Moreover, we have by  (ii) and (iii) of Lemma 2.2 with $ c_1:=\eta_2((1,\infty))$
$$(1-c) f_2(x):=c\sum_{n=1}^{\infty}\frac{c_1^ng_1^{n\otimes}(x)}{n!}\sim c_1g_1(x).$$ 
Thus, $f_2(x)$ is bounded and, by (i) of Lemma 2.2,
$ f_2(x) \in S_d$. 
Define $p(x)$ as
$$p(x):=cf_1(x)+(1-c)f_1\otimes f_2(x).$$ 
Then $p(x)$ is the density of $\mu(dx)=\xi_1*\xi_2(dx)$.
Note that
$$\lim_{x \to \infty}\frac{f_1(x)}{f_2(x)}=\lim_{x \to \infty}\frac{\exp({b_0x)} f_1(x)}{\exp({b_0x)} f_2(x)}=0.$$
Since $ f_2(x)$ is bounded and $ f_2(x) \in S_d$, for $0 \le u \le x$ and for sufficiently large $x >0$, there are $C>0$ and $b_1 >0$ such that
$$\frac{f_2(x-u)}{f_2(x)}\le C\exp(b_1 u).$$
By using Fatou's lemma, we have by $ f_2(x) \in S_d$
\begin{equation}
 \begin{split}
& \liminf_{x \to \infty}\frac{f_1\otimes f_2(x)}{f_2(x)}\\ \nonumber
&\ge \lim_{N \to \infty}\int_0^N\liminf_{x \to \infty}\frac{f_2(x-u)}{f_2(x)}f_1(u)du =1\\   
 \end{split}
 \end{equation}
and
\begin{equation}
 \begin{split}
& \limsup_{x \to \infty}\frac{f_1\otimes f_2(x)}{f_2(x)}\\ \nonumber
&\le \lim_{N \to \infty}\limsup_{x \to \infty}\int_0^N\frac{f_2(x-u)}{f_2(x)}f_1(u)du + \lim_{N \to \infty}\limsup_{x \to \infty}\int_N^{x}\frac{f_2(x-u)}{f_2(x)}f_1(u)du \\   
&\le  \lim_{N \to \infty}\int_0^N\limsup_{x \to \infty}\frac{f_2(x-u)}{f_2(x)}f_1(u)du+\lim_{N \to \infty}\int_N^{\infty} C\exp(b_1 u)f_1(u)du \\
&=1
\end{split}
 \end{equation} 
 Thus we have $p(x) \sim (1-c)f_1\otimes f_2(x)\sim (1-c)f_2(x)$ and hence, by (i) of Lemma 2.2, $p(x)\in S_d$, that is, (1).\hfill $\Box$

\medskip
 
Proof of Remark 1.1. Let $p_1(x)$ be the self-decomposable density  on $\mathbb R_+$ with 
L\'evy measure $1_{(0,1)}(x)x^{-1}k(x)dx$. By virtue of Theorem 25.3 of Sato \cite{s}, we have, for any $b>0$, 
$$\int_{0-}^{\infty}\exp(bx)p_1(x)dx<\infty.$$
Note from Remark 53.7 of Sato \cite{s} that $p_1(x)$ is continuous on $(0,\infty)$. Hence, we see from the unimodality of $p_1(x)dx$ by Theorem 53.1 of Sato \cite{s} that, for any  $b >0$,
$$\lim_{x \to \infty} \exp(bx)p_1(x)\le\lim_{x \to \infty}\int_{x-1}^{x}\exp(2bu)p_1(u)du=0.$$
 Let $\mu_2(dx):=c\delta_0(dx)+ (1-c)p_2(x)dx$ be the compound Poisson  distribution  on $\mathbb R_+$ with 
L\'evy measure $1_{(0,1)}(x)h(x)dx$. Define $c_1h_1(x):=1_{(0,1)}(x)h(x)$ with $c_1=\int_0^1h(x)dx$ Then $c:=\exp(-\int_0^1h(x)dx)$ and
$$(1-c)p_2(x):=c\sum_{n=1}^{\infty}\frac{c_1^nh_1^{n\otimes}(x)}{n!}.$$
Since $h(x)$ is bounded, $p_2(x)$ is bounded.
For $b >0$, we define the exponential tilt $(\mu_2)_{\langle b\rangle}$ of $\mu_2$ as
$$(\mu_2)_{\langle b\rangle}:=\frac{\exp(bx)}{\int_{0-}^{\infty}\exp(bx)\mu_2(dx)}\mu_2(dx).$$
Then, we find from Example 33.15 of Sato \cite{s} that $(\mu_2)_{\langle b\rangle}$ is the compound Poisson distribution with 
L\'evy measure 
$$\eta(dx):=1_{(0,1)}(x)\exp(bx)h(x)dx.$$ 
Note that the support of $\eta^{n*}$ is included in the interval $[0,n]$ and $\exp(bx)h(x)$ is bounded on $[0,1]$. Thus, with $C_0,C_1 >0$ we have, for any $b>0$,
$$\lim_{x \to \infty} \exp(bx)p_2(x) \le C_0 \lim_{x \to \infty}\sum_{n\ge x}\frac{C_1^n}{n!} =0$$
We have
$$  p_1\otimes p_2(2x)
\le\sup_{x\le u\le 2x}p_2(u)\int_0^x p_1(u)du+\sup_{x\le u\le 2x}p_1(u)\int_x^{2x} p_2(2x-u)du.$$
We find that
$$\lim_{x \to \infty}\exp(2bx)\sup_{x\le u\le 2x}p_2(u)\int_0^xp_1(u)du=0$$
and by the unimodality of $p_1(x)dx$
$$\lim_{x \to \infty}\exp(2bx)\sup_{x\le u\le 2x}p_1(u)\int_x^{2x}p_2(2x-u)du=0.$$
Thus we see that, for any  $b >0$,
$$\lim_{x \to \infty}\exp(2bx)p_1\otimes p_2(2x)=0$$
and hence, for any  $b >0$,
$$\lim_{x \to \infty}\exp(bx)f_1(x)=\lim_{x \to \infty}\exp(bx)(cp_1(x)+ (1-c)p_1\otimes p_2(x))= 0.$$
We have proved the remark. \hfill $\Box$

\medskip

Proof of Corollary 1.1. Assume that, for every $t >0$, there is some $b_0(t) >0$ such that
 $$\lim_{x \to \infty} \exp(b_0(t)x)f_1^t(x)=0$$ 
 and that $g_1(x)$ is bounded and $g_1(x) \in \mathcal{L}_{d}$. First, we prove assertion (i). Suppose that  $p^t(x) \in \mathcal{S}_d$ for some $t >0$. Then, we obtain from Theorem 1.1 that $g_1(x) \in \mathcal{S}_{d}$ and hence $p^t(x) \in \mathcal{S}_d$ for all $t >0$. Moreover, we find again from Theorem 1.1 that, for all $t >0$,
$$p^t(x)\sim t  C_0 g_1(x)\sim t p^1(x)$$
 with $C_0 :=\int_1^{\infty}g(x)dx$.
Next, we prove assertion (ii). Suppose that $p^1(x) \in \mathcal{L}_d$ and, for some 
$t \in  (0,1) \cup (1,\infty)$, there is 
$ C(t) \in (0,\infty)$            
 such that (1.2) holds. Then, we have $\mu \in \mathcal{L}_{loc}$ and (2.5) holds. Thus, we obtain from Lemma 2.5 that $C(t)=t$ and $\mu \in \mathcal{S}_{loc}$ and hence, by $p^1(x) \in \mathcal{L}_d$ and (ii) of Lemma 2.1, $p^1(x) \in \mathcal{S}_d$.\hfill $\Box$
 
 \medskip
Proof of Theorem 1.2.  Assume that $\int_{-\infty}^{\infty}|\widehat\mu(z)|dz < \infty$. We see from (ii) of Proposition 1.1 that (1), (3), and (4) are equivalent. Further, we find from (i) of Proposition 1.1 that (1) implies (2). Finally, we prove that (2) implies (1). Suppose that (2) holds. We see from Lemma 2.4 that $\mu_1 \in  \mathcal{S}_{loc} \subset \mathcal{L}_{loc}$ for every $c_1 > 0$. By (i) of Lemma 3.1, we can choose sufficiently large $c_1>0$ so that $\int_{-\infty}^{\infty}|\widehat\zeta_1(z)|dz < \infty$. Let $q(x)$ be the continuous density of $\zeta_1$. Then, we obtain from (ii) of Lemma 3.1 that, for every $\gamma > 0$,  $\lim_{x \to \infty}e^{\gamma x}q(x)=0.$ Since 
$$p(x)=\int_{0-}^{x+}q(x-u)\mu_1(du),$$
we see from (ii) of Lemma 2.3 that $p(x) \in \mathcal{L}_d$ and hence from (i) of Proposition 1.1 that $p(x) \in \mathcal{S}_d$.\hfill $\Box$

\medskip

Proof of Corollary 1.2. Assume that $\int_{-\infty}^{\infty}|\widehat\mu(z)|^t dz < \infty$ for every $t >0$. First, we prove assertion (i). Suppose that  $p^t(x) \in \mathcal{S}_d$ for some $t >0$. Then, we obtain from Theorem 1.2 that $\nu_{(1)} \in \mathcal{S}_{loc}$ and hence $p^t(x) \in \mathcal{S}_d$ for all $t >0$. Moreover, we find again from Theorem 1.2 that, for all $t >0$,
$$p^t(x)\sim t \nu((x,x+1]) \sim t p^1(x).$$
Next, we prove assertion (ii). Suppose that $p^1(x) \in \mathcal{L}_d$ and, for some 
$t \in  (0,1) \cup (1,\infty)$, there is 
$ C(t) \in (0,\infty)$            
 such that (1.3) holds. Then, we have $\mu \in \mathcal{L}_{loc}$ and (2.5) holds. Thus, we obtain from Lemma 2.5 that $C(t)=t$ and $\mu \in \mathcal{S}_{loc}$ and hence, by $p^1(x) \in \mathcal{L}_d$ and (ii) of Lemma 2.1, $p^1(x) \in \mathcal{S}_d$.\hfill $\Box$


\medskip
 

 
\medskip
\begin{thebibliography}{99}
\bibitem{afk} S. Asmussen, S. Foss, and D. Korshunov,  Asymptotics for sums of random variables with local subexponential behaviour, {\it J. Theoret. Probab.}, {\bf 16}, 489-518 (2003). 
\bibitem{cnw} J. Chover, P. Ney, and S. Wainger,  Functions of probability measures, {\it J. Analyse Math.}, {\bf 26},  255-302 (1973). 
\bibitem{egv} P. Embrechts, C.M. Goldie, and  N. Veraverbeke,  Subexponentiality and infinite divisibility, {\it  Z. Wahrscheinlichkeitstheorie Verw. Gebiet.}, {\bf 49},  335-347 (1979).
\bibitem{fkz} S. Foss, D. Korshunov, and  S. Zachary, {\it  An introduction to heavy-tailed and subexponential distributions}, Second edition. Springer Series in Operations Research and Financial Engineering. Springer, New York, (2013). 
\bibitem{jwcc} T. Jiang, Y. Wang, Z. Cui, and Y. Chen, On the almost decrease of a subexponential density, {\it  Statist. Probab. Letters}, {\bf 153},  71-79 (2019).
\bibitem{p} A.G. Pakes,  Convolution equivalence and infinite divisibility, {\it J. Appl. Probab.}, {\bf 41},  407-424 (2004).
 \bibitem{s} K. Sato,  {\it L\'evy processes and infinitely divisible distributions}, Cambridge Studies in Advanced Mathematics, 68 Cambridge Univ. Press. (2013). 
\bibitem{sw} T. Shimura and T. Watanabe, Infinite divisibility and generalized subexponentiality, {\it  Bernoulli}, {\bf 11},  445-469 (2005). 
\bibitem{w} T. Watanabe,  Convolution equivalence and distributions of random sums, {\it  Probab. Theory Related Fields}, {\bf 142}, (2008) 367-397. 
\bibitem{w3} T. Watanabe,  Asymptotic properties of Fourier transforms of $b$-decomposable distributions,  {\it J. Fourier Anal. Appl.}, {\bf 18},  803-827 (2012).

\bibitem{wy2}  T. Watanabe and K. Yamamuro, Local subexponentiality and \\
self-decomposability, {\it J. Theoret. Probab.}, {\bf 23},  1039-1067 (2010). 
\bibitem{wy3} T. Watanabe and K. Yamamuro, Tail behaviors of semi-stable distributions, {\it J. Math. Anal. Appl.,} {\bf 393},  108-121 (2012).

\bibitem{wy4} T. Watanabe and K. Yamamuro,  Two non-closure properties on the class of subexponential densities, {\it J. Theoret. Probab.}, {\bf 30},  1059-1075 (2017).
\end{thebibliography}
\end{document}